\documentclass[a4paper,12pt,final]{amsart}
\usepackage{times,a4wide,mathrsfs,amssymb,amsmath,amsthm,enumerate,xypic}

\newcommand{\C}{\mathbb{C}}
\newcommand{\ZZ}{\mathbb{Z}}

\newcommand{\QQ}{\mathbb{Q}}

\newcommand{\PP}{\mathbb{P}}

\newcommand{\OO}{\mathcal O}

\newcommand{\Sy}{\mathfrak S}

\newcommand{\XX}{\mathcal X}
\newcommand{\YY}{\mathcal Y}

\newcommand{\TT}{\mathcal T}

\newcommand{\Zz}{\mathcal Z}

\newcommand{\MM}{\mathcal M}

\newcommand{\FF}{\mathcal F}

\newcommand{\gr}{\hbox{Gr}}
\newcommand{\wt}{\widetilde}
\newcommand{\rom}{\romannumeral}

\DeclareMathOperator{\aut}{Aut}
\DeclareMathOperator{\bir}{Bir}
\DeclareMathOperator{\ide}{id}

\DeclareMathOperator{\ima}{Im}

\newtheorem{theorem}{Theorem}[section]
\newtheorem{claim}[theorem]{Claim}
\newtheorem{lemma}[theorem]{Lemma}

\newtheorem{corollary}[theorem]{Corollary}
\newtheorem{proposition}[theorem]{Proposition}

\newtheorem{remark}[theorem]{Remark}
\newtheorem{definition}[theorem]{Definition}
\newtheorem{convention}{Conventions}

\newtheorem{nonumbering}{Theorem}

\newtheorem{nonumberingc}{Corollary}

\newtheorem{nonumberingt}{Acknowledgements}

\begin{document}

\author[Robert Laterveer]
{Robert Laterveer}

\address{Institut de Recherche Math\'ematique Avanc\'ee,
CNRS -- Universit\'e 
de Strasbourg,\
7 Rue Ren\'e Des\-car\-tes, 67084 Strasbourg CEDEX,
FRANCE.}
\email{robert.laterveer@math.unistra.fr}

\title[Algebraic cycles and LLSS eightfolds]{Algebraic cycles and Lehn--Lehn--Sorger--van Straten eightfolds}

\begin{abstract} This article is about Lehn--Lehn--Sorger--van Straten eightfolds $Z$, and their anti-symplectic involution $\iota$. When $Z$ is birational to the Hilbert scheme of points on a K3 surface, we give an explicit formula for the action of $\iota$ on the Chow group of $0$-cycles of $Z$. The formula is in agreement with the Bloch--Beilinson conjectures, and has some non-trivial consequences for the Chow ring of the quotient.
%We also prove a statement about intersecting algebraic cycles on $Z$ with certain lagrangian subvarieties.
\end{abstract}

\thanks{Supported by ANR grant ANR-20-CE40-0023.}

\maketitle

\section{Introduction}

Given a smooth projective variety $X$ over $\C$, let $A^i(X):=CH^i(X)_{\QQ}$ denote the Chow groups of $X$ (i.e. the groups of codimension $i$ algebraic cycles on $X$ with $\QQ$-coefficients, modulo rational equivalence).   

In this article, we consider the Chow ring of Lehn--Lehn--Sorger--van Straten (=LLSS) eightfolds. These eightfolds (defined in terms of twisted cubic curves contained in a given cubic fourfold $Y$) form a $20$-dimensional locally complete family of
hyperk\"ahler varieties, deformation equivalent to the Hilbert scheme of points on a K3 surface \cite{LLSS}, \cite{AL}, \cite{CLehn}. One notable feature of the construction is that any LLSS eightfold admits an anti-symplectic involution $\iota$ \cite{MLe}.

Our main result gives an explicit formula for the action of $\iota$ on the Chow group of $0$-cycles, for certain LLSS eightfolds:

   \begin{nonumbering}[=Theorem \ref{main2}] Let $Z=Z(Y)$ be a LLSS eightfold, and let $\iota\in\aut(Z)$ be the anti-symplectic involution. Assume $Z$ is birational to a Hilbert scheme $Z^\prime=S^{[4]}$ where $S$ is a K3 surface. Then the action of $\iota$ on $A^8(Z^\prime)$ is given by
     \[ \begin{split} \iota_\ast [x,y,z,t]= [x,y,z,t] &-2\bigl( [x,y,z,o] + [x,y,t,o]+[x,z,t,o] + [y,z,t,o]\bigr)\\ &+4\bigl( [x,y,o,o]+[x,z,o,o]+\ldots+[z,t,o,o]\bigr)\\ &-8\bigl( [x,o,o,o]+[y,o,o,o]+[z,o,o,o]+[t,o,o,o]\bigr)\\&+16[o,o,o,o]\ .\\
     \end{split}\]
     (Here, $x,y,z,t\in S$ and $o\in A^2(S)$ denotes the distinguished $0$-cycle of \cite{BV}.)
   \end{nonumbering}

This theorem applies to a countably infinite number of divisors in the moduli space of LLSS eightfolds (the condition on $Z$ is equivalent to asking that the cubic fourfold $Y$ has discriminant of the form $(6n^2-6n+2)/a^2$ for some $n,a\in\ZZ$ \cite[Proposition 5.3]{LPZ}. In particular, Theorem \ref{main2} applies to Pfaffian cubics). 

The statement of Theorem \ref{main2} is similar to a statement for double EPW sextics proven in \cite[Theorem 2]{LV2}, and answers a question raised there \cite[Introduction]{LV2}.

The proof of Theorem \ref{main2} consists of a two-step argument. First, we prove that for {\em any\/} LLSS eightfold $Z$ (not necessarily birational to a Hilbert scheme), there is equality
  \begin{equation}\label{S3} \iota_\ast=-\ide\colon\ \ S^{mob}_3 A^8(Z)\cap A^8_{hom}(Z)\ \to\ A^8(Z)\ ,\end{equation}
  where $S^{mob}_\ast$ is related to Voisin's orbit filtration on zero-cycles \cite{V14}, cf. Definition \ref{defsmob} below. (As explained in Remark \ref{expect}, one expects that the piece $S^{mob}_3 A^8(Z)\cap A^8_{hom}(Z)$ is related, via the Bloch--Beilinson conjectures, to $H^{2,0}(Z)$. Since $\iota$ acts as $-\ide$ on $H^{2,0}(Z)$, equality (\ref{S3}) is in keeping with this expectation.) Next, we restrict attention to LLSS eightfolds birational to a Hilbert scheme $Z^\prime$. In this case, Vial \cite{V6} has constructed a multiplicative Chow--K\"unneth decomposition, giving a bigrading on the Chow ring $A^\ast_{(\ast)}(Z^\prime)$ such that the piece $S^{mob}_3 A^8(Z)\cap A^8_{hom}(Z^\prime)$ is exactly $A^8_{(2)}(Z^\prime)$. Exploiting the good formal properties of this bigrading, one can deduce Theorem \ref{main2} from equality (\ref{S3}).
  
  Theorem \ref{main2} admits a reformulation that does not mention LLSS eightfolds:
  
  \begin{nonumbering}[=Theorem \ref{main25}] Let $X=S^{[4]}$ be a Hilbert scheme, where 
   $S$ is a K3 surface of Picard number 1 and degree $d$ such that $d>2$ and $d= (6n^2-6n+2)/a^2$ for some $n,a\in\ZZ$. 
   Then $\bir(X)=\ZZ/2\ZZ$ and the non-trivial element $\iota\in\bir(X)$ acts on $A^8(X)$ as in Theorem \ref{main2}.
  \end{nonumbering}

As another consequence of Theorem \ref{main2}, we can show that the Chow ring of the quotient of the involution behaves somewhat like the Chow ring of Calabi--Yau varieties:

\begin{nonumberingc}[=Corollary \ref{cor}] Let $Z$ be a LLSS eightfold birational to a Hilbert scheme $Z^\prime=S^{[4]}$ where $S$ is a K3 surface. Let $X:=Z/\langle\iota\rangle$ be the quotient under the anti-symplectic involution $\iota\in\aut(Z)$.
Then
  \[ \begin{split} &\ima \Bigl(  A^2(X)^{\otimes 4}\to A^8(X)\Bigr)
                          = \ima\Bigl( A^3(X)\otimes A^2(X)\otimes A^2(X)\otimes A^1(X)\Bigr)
                          =\QQ[h^8]\ ,\\
                          \end{split}\]
        where $h\in A^1(X)$ is an ample divisor.
        \end{nonumberingc}
        
 Corollary \ref{cor} is expected to be true for {\em all\/} LLSS eightfolds, and not merely those birational to a Hilbert scheme. Since for a general LLSS eightfold we don't know how to construct a multiplicative Chow--K\"unneth decomposition, this expectation seems difficult to prove.

% Using similar arguments, we can also prove a statement about intersecting with certain lagrangian subvarieties:
% 
% \begin{nonumbering}[=Theorem \ref{main3}] Let $Z=Z(Y)$ be a LLSS eightfold, where $Y$ is a smooth cubic not containing a plane and such that the Fano variety of lines $F(Y)$ is birational to a Hilbert square of a K3 surface.
% Let $L\subset Z$ be either the lagrangian embedding of $Y$ in $Z$ \cite{LLSS}, or the lagrangian subvariety defined by twisted cubics contained in a hyperplane section of $Y$ \cite{SS}.
% One has equality
%   \[ L\cdot h^2\cdot A^2(Z)= \QQ[h^8]\ \ \ \hbox{in}\ A^8(Z)\ \]
%  (where $h\in A^1(Z)$ is a hyperplane section). 
%  \end{nonumbering}
%  
%  Again, this is in agreement with the general Bloch--Beilinson philosophy: the group $A^2_{hom}(Z)$ (of homologically trivial codimension 2 cycles) is expected to be related to $H^{2,0}(Z)$, which restricts to zero on the lagrangian subvariety $L$; as such, one expects that $L\cdot A^2_{hom}(Z)=0$.
%  
% Theorem \ref{main3} again applies to a countably infinite number of divisors in the moduli space of LLSS eightfolds (the condition on $Y$ is equivalent to asking that $Y$ has discriminant of the form $(2n^2+2n+2)/a^2$ for some $a,n\in\ZZ$ \cite[Theorem 2]{A}. In particular, Theorem \ref{main3} applies to Pfaffian cubics). 
% 
% In all likelihood, Theorem \ref{main3} is true for {\em all\/} LLSS eightfolds; unfortunately, in the absence of a multiplicative Chow--K\"unneth decomposition for $F(Y)$ (and for $Z$), this expectation seems difficult to prove.

 \vskip0.6cm

\begin{convention} In this article, the word {\sl variety\/} will refer to a reduced irreducible scheme of finite type over $\C$. A {\sl subvariety\/} is a (possibly reducible) reduced subscheme which is equidimensional. 

{\bf All Chow groups will be with rational coefficients}: we will denote by $A_j(X)$ the Chow group of $j$-dimensional cycles on $X$ with $\QQ$-coefficients; for $X$ smooth of dimension $n$ the notations $A_j(X)$ and $A^{n-j}(X)$ are used interchangeably. 

The notations $A^j_{hom}(X)$, $A^j_{AJ}(X)$ will be used to indicate the subgroups of homologically trivial, resp. Abel--Jacobi trivial cycles.
For a morphism $f\colon X\to Y$, we will write $\Gamma_f\in A_\ast(X\times Y)$ for the graph of $f$.
The contravariant category of Chow motives (i.e., pure motives with respect to rational equivalence as in \cite{Sc}, \cite{MNP}) will be denoted 
$\MM_{\rm rat}$.

We will write $H^j(X)$ to indicate singular cohomology $H^j(X,\QQ)$.

Given an automorphism $\sigma\in\aut(X)$, we will write $A^j(X)^\sigma$ (and $H^j(X)^\sigma$) for the subgroup of $A^j(X)$ (resp. $H^j(X)$) invariant under 
$\sigma$.
\end{convention}

\section{Preliminaries}

\subsection{Quotient varieties}
\label{ssquot}

\begin{definition} A {\em projective quotient variety\/} is a variety
  \[ X=Z/G\ ,\]
  where $Z$ is a smooth projective variety and $G\subset\hbox{Aut}(Z)$ is a finite group.
  \end{definition}
  
 \begin{proposition}[Fulton \cite{F}]\label{quot} Let $X$ be a projective quotient variety of dimension $n$. Let $A^\ast(X)$ denote the operational Chow cohomology ring. The natural map
   \[ A^i(X)\ \to\ A_{n-i}(X) \]
   is an isomorphism for all $i$.
   \end{proposition}
   
   \begin{proof} This is \cite[Example 17.4.10]{F}.
      \end{proof}

\begin{remark} It follows from Proposition \ref{quot} that the formalism of correspondences goes through unchanged for projective quotient varieties (this is also noted in \cite[Example 16.1.13]{F}). We may thus consider motives $(X,p,0)\in\MM_{\rm rat}$, where $X$ is a projective quotient variety and $p\in A^n(X\times X)$ is a projector. For a projective quotient variety $X=Z/G$, one readily proves (using Manin's identity principle) that there is an isomorphism
  \[  h(X)\cong h(Z)^G:=(Z,\Delta_G,0)\ \ \ \hbox{in}\ \MM_{\rm rat}\ ,\]
  where $\Delta_G$ denotes the idempotent ${1\over \vert G\vert}{\sum_{g\in G}}\Gamma_g$.  
  \end{remark}
  
    \subsection{LLSS eightfolds}

In this subsection we briefly recall the construction of the Lehn--Lehn--Sorger--van Straten eightfold (see \cite{LLSS} for additional details). Let $Y \subset \PP^5$ be a smooth cubic fourfold which does not contain a plane. Twisted cubic curves on $Y$ belong to an irreducible component $M_3(Y)$ of the Hilbert scheme $\textrm{Hilb}^{3m+1}(Y)$; in particular, $M_3(Y)$ is a ten-dimensional smooth projective variety and it is referred to as the Hilbert scheme of \emph{generalized twisted cubics} on $Y$. There exists a hyperk\"ahler eightfold $Z = Z(Y)$ and a morphism $u: M_3(Y) \to Z$ which factorizes as $u = \Phi \circ a$, where $a: M_3(Y) \to Z'$ is a $\PP^2$-bundle and $\Phi: Z' \to Z$ is the blow-up of the image of a Lagrangian embedding $j: Y\hookrightarrow Z$. By \cite{AL} (or alternatively \cite{CLehn}, or \cite[Section 5.4]{KLSV}), the manifold $Z$ is of $K3^{[4]}$-type.

For a point $j(y) \in j(Y) \subset Z$, a curve $C \in u^{-1}(j(y))$ is a singular cubic (with an embedded point in $y$), cut out on $Y$ by a plane tangent to $Y$ in $y$. In particular, $u^{-1}(j(Y)) \subset M_3(Y)$ is the locus of non-Cohen--Macaulay curves on $Y$. If instead we consider an element $C \in M_3(Y)$ such that $u(C) \notin j(Y)$, let $\Gamma_C := \langle C \rangle \cong \PP^3$ be the convex hull of the curve $C$ in $\PP^5$. Then, the fiber $u^{-1}(u(C))$ is one of the $72$ distinct linear systems of aCM twisted cubics on the smooth integral cubic surface $S_C := \Gamma_C \cap Y$ (cf.\ \cite{BK}). This two-dimensional family of curves, whose general element is smooth, is determined by the choice of a linear determinantal representation $[A]$ of the surface $S_C$ (see \cite[Chapter 9]{dolgachev}). This means that $A$ is a $3 \times 3$-matrix with entries in $H^0(\mathcal{O}_{\Gamma_C}(1))$ and such that $\det(A) = 0$ is an equation for $S_C$ in $\Gamma_C$. The orbit $[A]$ is taken with respect to the action of $\left(\textrm{GL}_3(\mathbb{C}) \times \textrm{GL}_3(\mathbb{
C})\right)/\Delta$, where $\Delta = \left\{ (tI_3,tI_3) \mid t \in \mathbb{C}^*\right\}$. Moreover, any curve $C' \in u^{-1}(u(C))$ is such that $I_{C'/S_C}$ is generated by the three minors of a $3 \times 2$-matrix $A_{0}$ of rank two, whose columns are in the $\mathbb{C}$-linear span of the columns of $A$. 

\subsection{Voisin's rational map and the involution $\iota$}
\label{ss:psi}
  
  \begin{proposition}[Voisin {\cite{V14}}]\label{psi} Let $Y\subset\PP^5$ be a smooth cubic fourfold not containing a plane. Let $F=F(Y)$ be the Fano variety of lines and let $Z=Z(Y)$ be the LLSS eightfold of $Y$. There exists a degree $6$ dominant rational map
  \[ \psi\colon\ \ F\times F\ \dashrightarrow\ Z\ .\]
%  This has the property that the composition
%   \[      H^{2,0}(F)\ \xrightarrow{ (pr_1)^\ast -(pr_2)^\ast}\ H^{2,0}(F\times F)\ \xrightarrow{\psi_\ast}\ H^{2,0}(Z) \]
%  is an isomorphism (where $pr_j$ denotes projection on the $j$--th factor). 

Moreover, let $\omega_F, \omega_Z$ denote the symplectic forms of $F$, resp.\ $Z$. Then we have
  \[ \psi^\ast(\omega_Z)= (p_1)^\ast(\omega_F) -(p_2)^\ast(\omega_F)\ \ \ \hbox{in}\ H^{2,0}(F\times F)\ \]
(where $p_j$ denotes projection on the $j$-th factor). 
   \end{proposition}

\begin{proof} This is \cite[Proposition 4.8]{V14}. Let us briefly recall the geometric construction of the rational map $\psi$. Let $(l,l') \in F \times F$ be a generic point, so that the two lines $l, l'$ on $Y$ span a three-dimensional linear space in $\mathbb{P}^5$. For any point $x \in l$, the plane $\langle x, l'\rangle$ intersects the smooth cubic surface $S = \langle l, l' \rangle \cap Y$ along the union of $l'$ and a residual conic $Q'_x$, which passes through $x$. Then, $\psi(l,l') \in Z$ is the point corresponding to the two-dimensional linear system of twisted cubics on $S$ linearly equivalent to the rational curve $l \cup_x Q'_x$ (this linear system actually contains the $\mathbb{P}^1$ of curves $\left\{ l \cup_x Q'_x \mid x \in l\right\}$).
\end{proof}

  \begin{remark} The indeterminacy locus of the rational map $\psi$ is the codimension $2$ locus $I\subset F\times F$ of intersecting lines (see \cite[Theorem 1.2]{Mura}).
  As proven in \cite[Theorem 1.1]{Chen}, the indeterminacy of the rational map $\psi$ is resolved by a blow-up with center $I\subset F\times F$ .
  \end{remark}

\begin{remark} \label{rmk: compatibility involution}
As is known from unpublished work of Lehn, Lehn, Sorger and van Straten, the LLSS eightfold $Z$ is equipped with a biholomorphic anti-symplectic involution $\iota \in \aut(Z)$ (cf. \cite{MLe}, \cite{Le}) which is induced, via the rational map $\psi: F \times F \dashrightarrow Z$, by the involution of $F \times F$ exchanging the two factors. Indeed, by \cite[Remark 3.8]{LLMS}, the involution $\iota$ acts on the subset $Z \setminus j(Y)$ as follows. A point $p \in Z \setminus j(Y)$ corresponds to a linear determinantal representation $[A]$ of a smooth cubic surface $S = Y \cap \mathbb{P}^3$, or equivalently to a \emph{six} of lines $(e_1, \ldots, e_6)$ on $S$ (see \cite[Section 9.1.2]{dolgachev} and \cite{BK}). Then, $\iota(p)$ is the point of $Z \setminus j(Y)$ associated with the (unique) six $(e'_1, \ldots, e'_6)$ on the same surface $S$ which forms a \emph{double-six} together with $(e_1, \ldots, e_6)$. In terms of linear determinantal representations, $\iota(p)$ corresponds to the transposed representation $[A^t]$ of the surface $S$ by \cite[Proposition 3.2]{BK}. It can be checked that, if $C$ is a twisted cubic curve on $S$ in the linear system parametrized by the point $p$, then, for any quadric $\mathcal{Q} \subset \langle C \rangle$ containing $C$, the residual intersection $C'$ of $S$ with $\mathcal{Q}$ belongs to the fiber $u^{-1}(\iota(p))$ (see for instance \cite[Footnote 47]{debarre}). By the geometric description of the rational map $\psi$, we conclude that $\iota(\psi(l,l')) = \psi(l',l)$: in fact, for $x \in l$ and $x' \in l'$, the reducible quadric $\mathcal{Q} = \langle x, l'\rangle \cup \langle x', l \rangle$ intersects the surface $S$ along the union of $l \cup_x Q'_x$ and $l' \cup_{x'} Q_{x'}$.
\end{remark}

 \subsection{Voisin's orbit filtration}
 \label{ss:vois}

 \begin{definition}[\cite{V14}]\label{defs}
 Let $X$ be a variety of dimension $n$. One defines $S_j X\subset X$ as the set of points whose orbit under rational equivalence has dimension $\ge j$.
 The descending filtration $S_\ast$ on $A^n(X)$ is defined by letting $S_j A^n(X)$ be the subgroup generated by points $x\in S_j X$.
 \end{definition}
 
 \begin{remark} Let $X$ be a hyperk\"ahler variety. As explained in \cite{V14}, the expectation is that the filtration $S_\ast$ is opposite to the (conjectural) Bloch--Beilinson filtration, and thus provides a splitting of the Bloch--Beilinson filtration on $A^n(X)$ (this is motivated by Beauville's famous speculations in \cite{Beau3}). This splitting should be ``natural'', in the sense that it can be defined by correspondences and that it induces a bigrading on the Chow ring of $X$.
 \end{remark}

\begin{definition}\label{defsmob} Let $Z$ be a LLSS eightfold, and let $D\subset Z$ be the uniruled divisor obtained by resolving the indeterminacy locus of the rational map $\psi$ of Proposition \ref{psi}. One defines $S_j^{mob} A^8(Z)$ to be the subgroup generated by points in $S_j (Z\setminus D)$.
\end{definition}

\begin{remark}\label{expect} There is an inclusion 
  \[ S_j^{mob} A^8(Z)\subset S_j A^8(Z)\ ,\] 
which possibly is an equality (this is true when $Z$ is birational to a Hilbert scheme of points on a K3 surface, cf. the proof of Lemma \ref{incl} below). As explained in \cite{V14}, the piece $S_3 A^8(Z)\cap A^8_{hom}(Z)$ is supposed to be isomorphic to the graded $\gr_{F_{BB}}^{2} A^8(Z)$, where $F^\ast_{BB}$ is the (conjectural) Bloch--Beilinson filtration. On the other hand, the graded $\gr_{F_{BB}}^{2} A^8(Z)$ should be defined by the projector $\pi^{14}_Z$ of some Chow--K\"unneth decomposition for $Z$.
In particular, the expectation is that $S_3^{mob} A^8(Z)$ is somehow related to ($H^{14}(Z)$ and hence) $H^{2,0}(Z)$.
\end{remark}

\subsection{MCK decomposition}

\begin{definition}[Murre \cite{Mur}] Let $X$ be a smooth projective variety of dimension $n$. We say that $X$ has a {\em CK decomposition\/} if there exists a decomposition of the diagonal
   \[ \Delta_X= \pi_0+ \pi_1+\cdots +\pi_{2n}\ \ \ \hbox{in}\ A^n(X\times X)\ ,\]
  such that the $\pi_i$ are mutually orthogonal idempotents and $(\pi_i)_\ast H^\ast(X)= H^i(X)$.
  
  (NB: ``CK decomposition'' is shorthand for ``Chow--K\"unneth decomposition''.)
\end{definition}

\begin{remark} The existence of a CK decomposition for any smooth projective variety is part of Murre's conjectures \cite{Mur}, \cite{J4}. 
%If a quotient variety $X$
%has finite--dimensional motive, and the K\"unneth components are algebraic, then $X$ has a CK decomposition (this can be proven just as \cite{J2}, where this is stated for smooth $X$).
\end{remark}

\begin{definition}[Shen--Vial \cite{SV}] Let $X$ be a smooth projective variety of dimension $n$. Let $\Delta_X^{sm}\in A^{2n}(X\times X\times X)$ be the class of the small diagonal
  \[ \Delta_X^{sm}:=\bigl\{ (x,x,x)\ \vert\ x\in X\bigr\}\ \subset\ X\times X\times X\ .\]
  An {\em MCK decomposition\/} is a CK decomposition $\{\pi^X_i\}$ of $X$ that is {\em multiplicative\/}, i.e. it satisfies
  \[ \pi^X_k\circ \Delta_X^{sm}\circ (\pi^X_i\times \pi^X_j)=0\ \ \ \hbox{in}\ A^{2n}(X\times X\times X)\ \ \ \hbox{for\ all\ }i+j\not=k\ .\]
  
 (NB: ``MCK decomposition'' is shorthand for ``multiplicative Chow--K\"unneth decomposition''.) 
  
 A {\em weak MCK decomposition\/} is a CK decomposition $\{\pi^X_i\}$ of $X$ that satisfies
    \[ \Bigl(\pi^X_k\circ \Delta_X^{sm}\circ (\pi^X_i\times \pi^X_j)\Bigr){}_\ast (a\times b)=0 \ \ \ \hbox{for\ all\ } a,b\in\ A^\ast(X)\ .\]
  \end{definition}
  
  \begin{remark} The small diagonal (seen as a correspondence from $X\times X$ to $X$) induces the {\em multiplication morphism\/}
    \[ \Delta_X^{sm}\colon\ \  h(X)\otimes h(X)\ \to\ h(X)\ \ \ \hbox{in}\ \MM_{\rm rat}\ .\]
 Suppose $X$ has a CK decomposition
  \[ h(X)=\bigoplus_{i=0}^{2n} h^i(X)\ \ \ \hbox{in}\ \MM_{\rm rat}\ .\]
  By definition, this decomposition is multiplicative if for any $i,j$ the composition
  \[ h^i(X)\otimes h^j(X)\ \to\ h(X)\otimes h(X)\ \xrightarrow{\Delta_X^{sm}}\ h(X)\ \ \ \hbox{in}\ \MM_{\rm rat}\]
  factors through $h^{i+j}(X)$.
  
  If $X$ has a weak MCK decomposition, then setting
    \[ A^i_{(j)}(X):= (\pi^X_{2i-j})_\ast A^i(X) \ ,\]
    one obtains a bigraded ring structure on the Chow ring: that is, the intersection product sends $A^i_{(j)}(X)\otimes A^{i^\prime}_{(j^\prime)}(X) $ to  $A^{i+i^\prime}_{(j+j^\prime)}(X)$.
    
      It is expected that for any $X$ with a weak MCK decomposition, one has
    \[ A^i_{(j)}(X)\stackrel{??}{=}0\ \ \ \hbox{for}\ j<0\ ,\ \ \ A^i_{(0)}(X)\cap A^i_{hom}(X)\stackrel{??}{=}0\ ;\]
    this is related to Murre's conjectures B and D, that have been formulated for any CK decomposition \cite{Mur}.

  The property of having a MCK decomposition is severely restrictive, and is closely related to Beauville's ``(weak) splitting property'' \cite{Beau3}. For more ample discussion, and examples of varieties with a MCK decomposition, we refer to \cite[Section 8]{SV}, as well as \cite{V6}, \cite{SV2}, \cite{FTV}, \cite{LV}.
    \end{remark}

In what follows, we will make use of the following:

\begin{theorem}[Vial \cite{V6}]\label{K3m} Let $X$ be a Hilbert scheme $S^{[m]}$, where $S$ is a $K3$ surface. Then $X$ has a MCK decomposition. Let $A^\ast_{(\ast)}(X)$ denote the resulting bigrading. One has
  \[  \bigoplus_{j=0}^r A^{2m}_{(2j)}(X) = S_{m-r} A^{2m}(X)\ ,\]
  where $S_\ast$ is Voisin's orbit filtration (Definition \ref{defs}).
  
 The pieces $A^{2m}_{(\ast)}(X)$ can be described as follows: $A^{2m}_{(0)}(X)$ is generated by $[o,o,\ldots,o]$ (where $o\in A^2(S)$ is the distinguished $0$-cycle of \cite{BV}).  $A^{2m}_{(2)}(X)$ is generated by expressions 
   \[ [x,o,\ldots,o]-[y,o,\ldots,o] \] 
   (where $x,y$ range over all points of $S$). $A^{2m}_{(4)}(X)$ is generated by expressions 
  \[ [x,y,o,\ldots,o]-[x,o,\ldots,o] -[y,o,\ldots,o] +[o,o,\ldots,o]\ .\]
  $A^{2m}_{(6)}(X)$ is generated by expressions 
  \[ \begin{split}   [x,y,z,o,\ldots,o]-[x,y,o,\ldots,o] -[x,z,o,\ldots,o] -[y,z,o,\ldots,o]  +&\\ 
 [ x,o,\ldots,o] + [y,o,\ldots,o] + [z,o,\ldots,o]   - [o,o,\ldots,o]\ .&\\\ 
 \end{split}\]  
  (And so on...)
    \end{theorem}

\begin{proof} The MCK decomposition is constructed in \cite[Theorem 1]{V6} (an alternative construction is given in \cite{NOY}). The compatibility with Voisin's filtration $S_\ast$ is stated in \cite[End of Section 4.1]{V14}. The formulae for the generators follow from \cite[Theorem 2.5]{V14}, where it is proven that the filtration $S_\ast$ coincides with the filtration $N_\ast$ defined in \cite[Proposition 2.2]{V14}.

(The $m=2$ case also follows from \cite{SV}.)
\end{proof}

\subsection{Multiplicative behaviour of the Chow ring of $S^{[m]}$}

Let $X$ be a Hilbert scheme of points on a K3 surface. The following proposition shows that the Chow group of $0$-cycles on $X$ is completely determined by the piece $A^2_{(2)}(X)$. (This can be interpreted as a motivic manifestation of the fact that $H^\ast(X,\OO_X)$ is determined by $H^2(X,\OO_X)$, since $A^2_{(2)}(X)$ is defined by the projector $\pi^{2}_X$).

\begin{proposition}\label{mult} Let $X$ be a Hilbert scheme $X=S^{[m]}$ where $S$ is a K3 surface, and let $A^\ast_{(\ast)}(X)$ refer to the bigrading of Theorem \ref{K3m}.

\item{(\rom1)} For any $j\ge 1$, intersection product induces surjections
  \[  A^2_{(2)}(X)^{\otimes j}\ \twoheadrightarrow\ A^{2j}_{(2j)}(X)\ .\]
  
  \item{(\rom2)} Let $h\in A^1(X)$ be a big divisor. For any $j\ge 0$ there are isomorphisms
  \[ \cdot h^{2m-j}\colon\ \ A^j_{(j)}(X)\ \xrightarrow{\cong}\ A^{2m}_{(j)}(X)\ .\]
  \end{proposition}
  
  \begin{proof}
  
  \noindent
  (\rom1) Let us just treat the case $j=2$ (the other cases are only notationally more complicated). 
  Let us write $A^\ast(S^m)^{\Sy_m}$ for the $\Sy_m$-invariant part, where the symmetric group on $m$ elements $\Sy_m$ acts on $S^m$ by permutation of the factors.
    We need to show that the map induced by intersection product
    \[ \Bigl(A^2_{(2)}(S^m)^{\Sy_m}\Bigr)^{\otimes 2}\ \to\ A^4_{(4)}(S^m)^{\Sy_m}\  \]
    is surjective. By definition, an element $a$ on the right-hand side is a symmetric cycle of the form
    \[ \bigl( \pi^2_S\times \pi^2_S\times \pi^0_S\times\cdots\times \pi^0_S + \hbox{permutations}\bigr){}_\ast (d)\ \]
    (where $d\in A^4(S^m)$). That is, 
    \[ a = \sum_{\sigma\in\Sy_m} \sigma_\ast (a_1\times a_2\times S\times S\times \cdots\times S)  \ \ \hbox{in}\ A^4(S^m)\ ,\]
    where $a_1, a_2\in A^2(S)$ are degree zero $0$-cycles.
    
    Let us define
    \[ \begin{split}  b_1&:= \sum_{\sigma\in\Sy_m} \sigma_\ast (a_1\times S\times S\times \cdots\times S)  \ \ \in\ A^2_{(2)}(S^m)^{\Sy_m}\ ,\\
                           b_2&:= \sum_{\sigma\in\Sy_m} \sigma_\ast (S\times a_2\times S\times \cdots\times S)  \ \ \in\ A^2_{(2)}(S^m)^{\Sy_m}\ .\\
                          \end{split}\]
        We have that
        \[ b_1\cdot b_2= \sum_{\sigma\in\Sy_m} \sum_{\sigma^\prime\in\Sy_m} (\sigma\circ\sigma^\prime)_\ast (a_1\times a_2\times S\times S\times \cdots\times S) \ \]
        is a non-zero multiple of $a$, and so we have proven surjectivity.
 
 \noindent       
  (\rom2) Again, one reduces to the cartesian product $S^m$. That is, one considers
    \[ \begin{array}[c]{ccc}
                      &&X\\
                        && \ \  \downarrow{\scriptstyle f}\\
                       S^m&\xrightarrow{p}&\ \ \ \ \ \ \ \ \ \ \ \ \ \ \ S^{(m)}=S^m/\Sy_m\\
                       \end{array}\]
                       where $f$ is the Hilbert--Chow morphism and $p$ is the quotient morphism.
                       
     There is a commutative diagram
     \[  \begin{array}[c]{ccc}
                  A^j_{(j)}(X) &\xrightarrow{\cdot h^{2m-j}}&  A^{2m}_{(j)}(X)\\
                  &&\\
                \ \   \downarrow {\scriptstyle f_\ast} &&  \ \   \downarrow {\scriptstyle f_\ast}\\
                &&\\
                A^j(S^{(m)}) &\xrightarrow{\cdot f_\ast (h)^{2m-j}}&  A^{2m}_{}(S^{(m)})\\
                &&\\
                     \ \   \downarrow {\scriptstyle p^\ast} &&  \ \   \downarrow {\scriptstyle p^\ast}\\
                     &&\\
                     A^j_{(j)}(S^m)^{\Sy_m} &\xrightarrow{\cdot  p^\ast f_\ast (h)^{2m-j}}&  A^{2m}_{(j)}(S^m)^{\Sy_m}\\
                     \end{array} \]
                     in which composition of two vertical arrows on the left resp. on the right is an isomorphism (NB: the commutativity of the upper square follows from the projection formula).
  
     Thus, it will suffice to prove that for a big $\Sy_m$-invariant divisor $h$ on $S^m$, there are isomorphisms
   \[      \cdot h^{2m-j}\colon\ \ A^j_{(j)}(S^m)^{\Sy_m}\ \xrightarrow{\cong}\ A^{2m}_{(j)}(S^m)^{\Sy_m}\ .\]
   The case $j=2$ is proven in \cite[Theorem 3.1]{epw}; the general case is only notationally more complicated.
     \end{proof}

\begin{remark} The $m=2$ case of Proposition \ref{mult} is already proven in \cite{SV}.
\end{remark}

 \subsection{Birational hyperk\"ahler varieties}
 
 \begin{theorem}[Rie\ss~\cite{Rie}]\label{rie} Let $\phi\colon X\dashrightarrow
		X^\prime$ be a birational map between two hyperK\"ahler varieties of dimension
		$n$. 
		\begin{enumerate}[(i)]
			\item  There exists a correspondence $R_\phi\in A^n(X\times X^\prime)$ such
			that 
			\[ (R_\phi)_\ast \colon \ \ A^\ast(X) \ \to\  A^\ast(X^\prime)\]
			is a graded ring isomorphism.
			\item For any $j$, there is equality
			\[ (R_\phi)_\ast c_j(X)= c_j(X^\prime)\ \  \hbox{ in}\  A^j(X^\prime)\ .\]
					\item
					 There is equality
		    \[ (\Gamma_R)_\ast=(\bar{\Gamma}_\phi)_\ast\colon\ \  
		    A^j(X)\ \to\ \ 
		 A^j(X^\prime)\ \ \ \hbox{for}\ j=0,1, n-1, n\ \]
		    (where $\bar{\Gamma}_\phi$ denotes the closure of the graph of $\phi$). 		
					%There exists a birational map $\phi^\prime\colon X\dashrightarrow X^\prime$ such that there is equality
			%$(R_\phi)_\ast=(\bar{\Gamma}_{\phi^\prime})_\ast\colon \A^n(X)\to \A^n(X^\prime)$,
			%$j=0,1, n-1, n$,
			%where $\bar{\Gamma}_{\phi^\prime}$ denotes the closure of the graph of $\phi^\prime$.  
			%\Charles{Je ne comprends pas\,: Que se passe-t-il si on prend $X=X'$, $\Gamma_R
			%= \mathrm{id}_X$ et $\phi$ un automorphisme non-trivial ?
			%Dans notre cas, je pense qu'on peut prendre des familles au-dessus de $C$ o\`u
			%l'involution est definie au-dessus de $C\backslash 0$, et on prend $\Gamma_R$
			%la specialisation a $0$ de l'involution.
			%
			%(Robert) J'ai modifie l'enonce, afin de bien faire ressortir (dans $(\rom1)$) que
			%la construction associe une correspondence a une application birationnelle (et
			%non pas a une paire de varieties birationnelles - observation elementaire sans
			%doute, mais qui n'etait pas clair dans ma tete jusqu'ici). Cela fait
			%disparaitre ton objection (car a un automorphisme $\psi\in\aut(X)$ la
			%construction associe $R_\psi=\Gamma_\psi$).
			%
			%(Charles) Bien vu, je n'avais jamais realis\'e que la correspondence est
			%associ\'ee \`a une application birationnelle.  }  
		\end{enumerate}
		%\noindent{(\rom1)}	 
		%	 There exists a correspondence $\Gamma_R\in \A^n(X\times X^\prime)$ such that 
		%    \[ (\Gamma_R)_\ast \colon\ \ \A^\ast(X)\ \to\ \A^\ast(X^\prime) \]
		%    is a graded ring isomorphism.
		%    
		%  \noindent{(\rom2)} For any $j$, there is equality
		%   \[ (\Gamma_R)_\ast c_j(X)= c_j(X^\prime)  \ \ \ \hbox{in}\ \A^j(X^\prime)\
		%.\]
		%    
		%   \noindent{(\rom3)}
		%   There is equality
		%    \[ (\Gamma_R)_\ast=(\bar{\Gamma}_\phi)_\ast\colon\ \ \A^j(X)\ \to\
		%\A^j(X^\prime)\ \ \ \hbox{for}\ j=0,1, n-1, n\ \]
		%    (where $\bar{\Gamma}_\phi$ denotes the closure of the graph of $\phi$).    
	\end{theorem}
	
	\begin{proof} Item $(\rom1)$ is \cite[Theorem 3.2]{Rie}, while item $(\rom2)$ is
		\cite[Lemma 4.4]{Rie}.
		
		In a nutshell, the construction of the correspondence $R_\phi$ is as follows: \cite[Section 2]{Rie} provides a diagram 
		$$\xymatrix{\XX \ar[dr] \ar@{-->}[rr]^{\Psi} & & \XX' \ar[dl] \\
			& C
		}$$
	where $\XX, \XX^\prime$ are algebraic spaces over a quasi-projective
	curve $C$ such that the fibers $\XX_0, \XX^\prime_0$ are isomorphic to $X$
	resp. $X^\prime$, and where $\Psi\colon \XX
	\dashrightarrow \XX^\prime$ is a birational map inducing an isomorphism $\XX_{C\setminus 0}\cong  
	\XX^\prime_{C\setminus 0}$
	 and whose  restriction $\Psi\vert_{\XX_0}$
	coincides with $\phi$.	
	The correspondence $R_\phi$ is then defined as the specialization (in the sense
	of \cite{F}, extended to algebraic spaces) of the graphs of the isomorphisms
	$\XX_c\cong\XX^\prime_c$, $c\not=0$.

Point $(\rom3)$ is not stated explicitly in \cite{Rie}, and can be seen as follows.
Let $U\subset X$ be the locus on which $\phi$ induces an isomorphism. The
		complement $T:=X\setminus U$ has codimension $\ge 2$. Using complete
		intersections of hypersurfaces, one can find a closed subset $\TT\subset\XX$ of
		codimension $2$ such that $\TT_0$ contains $T$, i.e. $\Psi$ restricts to an
		isomorphism on $V:=\XX_0\setminus \TT_0$.  
		Since specialization commutes with pullback, the restriction of $R_\phi$ to
		$V\times X^\prime$ is the specialization of the graphs of the morphisms
		$\XX_c\setminus\TT_c\to \XX^\prime_c$, $c\not=0$ to $V\times X^\prime$, which
		is exactly the graph of the morphism $\phi\vert_V$, i.e.
		\[ R_\phi\vert_{V\times X\prime}= \Gamma_{(\phi\vert_V)}= \bar{\Gamma}_\phi
		\vert_{V\times X^\prime}\ \ \hbox{in}\ A^n(V\times X^\prime)\ .\]
		It follows that one has
		\[ R_\phi = \bar{\Gamma}_\phi + \gamma\ \ \ \hbox{in}\  A^n(X\times
		X^\prime)\ ,\]
		where $\gamma$ is some cycle supported on $\TT_0\times X^\prime$. This proves
		point $(\rom3)$: $0$-cycles and $1$-cycles on $X$ can be moved to be
		disjoint of $\TT_0$, and so $\gamma_\ast$ acts as zero on $A^j(X)$ for $j\ge
		n-1$. Likewise, $\gamma^\ast$ acts as zero on $A^j(X^\prime)$ for $j\le 1$ for
		dimension reasons, and so (by inverting the roles of $X$ and $X^\prime$)
		statement $(\rom3)$ is proven. 
	\end{proof}

	\begin{remark}\label{birat} As observed in \cite[Introduction]{V6}, the
		correspondence $R_\phi$ of Theorem \ref{rie} actually induces an isomorphism of
		Chow motives as $\QQ$-algebra objects. This implies that the property ``having
		a MCK decomposition'' is birationally invariant among hyperk\"ahler varieties.
	\end{remark}

  \begin{remark}\label{birat2} Let $Z$ be a hyperk\"ahler variety birational to a Hilbert scheme $Z^\prime=S^{[m]}$. It follows from Theorem \ref{rie} that 
  Proposition \ref{mult} is valid for $Z$.
  \end{remark}

 \subsection{The Franchetta property} 
 
 An important role in this paper is played by the Franchetta property. In a nutshell, a family of varieties has the Franchetta property if algebraic cycles that exist for all members
 of the family are well-behaved. This is relevant here, because the graph of the anti-symplectic involution of a LLSS eightfold is such a cycle: it exists for all members of the locally complete family of LLSS eightfolds.

   \begin{definition} Let $\XX\to B$ be some family of smooth projective varieties $X_b$, $b\in B$. We say that $\XX\to B$ {\em has the Franchetta property\/} if for any $a\in A^i( \XX)$ and for any $b\in B$, there is equivalence
   \[ a\vert_{X_b}=0\ \ \hbox{in}\ H^{2i}(X_b)\ \ \iff\ \ a\vert_{X_b}=0\ \ \hbox{in}\ A^{i}(X_b)\  .\] 
   \end{definition}
   
  This property is studied in \cite{OG}, \cite{PSY} for the universal family of $K3$ surfaces of low genus. This is extended to certain families of hyperk\"ahler varieties in \cite{FLV}, \cite{FLV3}.

  \begin{proposition}\label{fpcubic} Let $\YY\to B$ denote the universal family of smooth cubic fourfolds, where $B$ is a Zariski open in $\PP^{55}$. The Franchetta property holds for
    \[ A^4(\YY\times_B \YY) \ .\]
    \end{proposition}
    
    \begin{proof} Let $\Gamma\in A^4(\YY\times_B \YY)$ be a correspondence that is homologically trivial when restricted to the very general fiber.
    According to \cite[Proposition 1.6]{V1} (cf. also \cite[Proposition 5.1]{LNP}, where the precise form of Voisin's argument applied here is detailed), there exists a cycle $\gamma\in A^4(\PP^5\times\PP^5)$ such that
     \[ \Gamma\vert_b = \gamma\vert_b\ \ \ \hbox{in}\ A^4(Y_b\times Y_b)\ \ \ \forall b\in B\ .\]
     The cycle $\gamma$ is of the form $\sum_{j=0}^4 c_j\, h^j\times h^{4-j}$, where $h\in A^1(\PP^5)$ is a hyperplane section and $c_j\in\QQ$.
     The assumptions imply that $\gamma\vert_b$ is homologically trivial for $b\in B$ very general. Since the summands $(h^j\times h^{4-j})\vert_b$ live in different summands of the K\"unneth decomposition of $H^8(Y_b\times Y_b)$, each $c_j\, (h^j\times h^{4-j})\vert_b$ must be homologically trivial, which is only possible when $c_j=0$, and so $\gamma=0$.     
     \end{proof} 
   
 \begin{remark} Using the Franchetta property for $\FF\times_B \FF$ (where $\FF$ is the universal Fano variety of lines in a cubic fourfold) \cite{FLV}, it is possible to prove that $\YY\times_B \YY\to B$ (and also $\YY^{(2)}\times_B \YY^{(2)}\to B$) has the Franchetta property in all codimensions. 
 We will not need this generality here.
 \end{remark}

 \subsection{The transcendental motive of a cubic fourfold}
   
   \begin{proposition}[\cite{BP}]\label{bp} Let $Y\subset\PP^5$ be a smooth cubic fourfold. There exists a Chow motive $t(Y)\in\MM_{\rm rat}$, with the property that
     \[ \begin{split}  &A^\ast(t(Y))= A^3_{hom}(Y)\ ,\\
                              &H^\ast(t(Y))= H^4_{tr}(Y)\ \\
                              \end{split} \]
            (where $H^4_{tr}(Y)$ denotes the orthogonal complement of the algebraic classes with respect to cup-product).
       \end{proposition}
       
       \begin{proof} The construction, which is a variant of the construction of the transcendental motive of a surface \cite{KMP}, is done in \cite[Section 2]{BP}. Alternatively, one can use the projectors on the algebraic part of cohomology constructed for any variety satisfying the standard conjectures in \cite{V2}.
        \end{proof}

  \section{Main results}
  
For any LLSS eightfold $Z$, we construct (Theorem \ref{motive}) a Chow motive $t(Z)$ which is responsible for one piece of Voisin's orbit filtration on the Chow group of zero-cycles, i.e. the motive $t(Z)$ determines $S_3^{mob} A^8(Z)\cap A^8_{hom}(Z)$.  
 Using this construction, we can show (Theorem \ref{main-1}) that the anti-symplectic involution $\iota$ acts as $-\ide$ on this piece $S_3^{mob} A^8(Z)\cap A^8_{hom}(Z)$. Next, we restrict to LLSS eightfolds birational to Hilbert schemes, where we have a MCK decomposition. Using the good properties of the MCK decomposition, we deduce (Theorems \ref{main0} and \ref{main2}) that the involution $\iota$ acts on zero-cycles in the expected way.

  \subsection{An isomorphism of motives}
  
  \begin{theorem}\label{motive} Let $Z=Z(Y)$ be a LLSS eightfold, where $Y\subset\PP^5$ is any smooth cubic not containing a plane.
  There exists a motive $t(Z)\in\MM_{\rm rat}$ with the property that
    \[ A^\ast(t(Z))=S_3^{mob} A^8(Z)\cap A^8_{hom}(Z)\ .\]
    There is an isomorphism of motives
    \[ t(Z)\cong t(Y)(-5)\ \ \ \hbox{in}\ \MM_{\rm rat}\ ,\]
  inducing in particular an isomorphism of $\QQ$-vector spaces
   \[   S_3^{mob} A^8(Z)\cap A^8_{hom}(Z)\cong A^3_{hom}(Y)\ .\]  
     \end{theorem}  
  
 \begin{proof} Let us write $F:=F(Y)$ for the Fano variety of lines, and $P\in A^3(F\times Y)$ for the universal line. We
 define a correspondence
   \begin{equation}\label{defQ} Q:= \bar{\Gamma}_\psi\circ ({}^t \Gamma_{p_1}-{}^t \Gamma_{p_1})\circ {}^t P\ \ \ \in\ A^3(Y\times Z)\ ,\end{equation}
   where $\psi\colon F\times F\dashrightarrow Z$ is Voisin's rational map, and $p_j\colon F\times F\to F$ is projection to the $j$th factor.
   The correspondence $Q$ has the following property:
   
   \begin{claim}\label{Q} Let $h\in A^1(Z)$ be the polarization. There exists a constant $c_h\in\QQ^\ast$ such that
     \[  c_h\, {}^t Q\circ \Gamma_{h^6}\circ Q\circ \pi^{tr}_Y=\pi^{tr}_Y\ \ \ \hbox{in}\ A^4(Y\times Y)\ .\]
     (Here $\Gamma_{h^6}:=\Gamma_\tau\circ {}^t \Gamma_\tau\in A^{14}(Z\times Z)$, where $\tau\colon {S}\to Z$ is the inclusion of a smooth complete intersection surface $S\subset Z$ cut out by $h$.)
     
     In particular, the composition
     \[ A^3_{hom}(Y)\ \xrightarrow{Q^\ast}\ A^2_{hom}(Z)\ \xrightarrow{\cdot h^6}\ A^8_{hom}(Z)\ \xrightarrow{(c_h Q)_\ast}\ A^3_{hom}(Y) \]
     is the identity.
        \end{claim}
 
 Admitting the claim, the theorem is readily proven: one defines the motive $t(Z)$ by the projector
   \[ \pi^{tr}_Z:=c_h\,  \Gamma_{h^6}\circ Q\circ \pi^{tr}_Y\circ  {}^t Q\ \ \ \in\ A^8(Z\times Z)\ .\]
   The claim implies that $\pi^{tr}_Z$ is idempotent, and that
   \[   {}^t Q\colon\ \    t(Z):=(Z,\pi^{tr}_Z,5)\ \to\ t(Y):=(Y,\pi^{tr}_Y,0)\ \ \ \hbox{in}\ \MM_{\rm rat} \]
   is an isomorphism, with inverse
   \[    c_h\, \Gamma_{h^6}\circ Q  \ \ \ \in\ A^{9}(Y\times Z)\ .\] 
 
 The isomorphism of motives $t(Z)\cong t(Y)(-5)$ implies that $A^\ast(t(Z))=A^3(t(Z))$. To determine $A^3(t(Z))$, we use the relation with the Fano variety of lines $F:=F(Y)$ via Voisin's rational map $\psi$. Indeed, by construction we have
   \[   (\pi^{tr}_Z)_\ast A^8(Z)= \ima\Bigl(  Q_\ast A^3_{hom}(Y)\ \xrightarrow{\cdot h^6}\ A^8(Z)\Bigr)\ .\]
  Now, there are equalities
      \[ \begin{split} &\ima\Bigl(  Q_\ast A^3_{hom}(Y)\ \xrightarrow{\cdot h^6}\ A^8(Z)\Bigr) \\
                       =&      \psi_\ast    \ima\Bigl(  \bigl((p_1)^\ast-(p_2)^\ast \bigr) P^\ast A^3_{hom}(Y)\ \xrightarrow{\cdot\psi^\ast(h^6) }\ A^8(F\times F)\Bigr)\\
                       =&      \psi_\ast   \ima \Bigl(   \bigl( A^2_{(2)}(F)\otimes A^0(F)\oplus A^0(F)\otimes A^2_{(2)}(F)\bigr)\ \xrightarrow{\cdot\psi^\ast(h^6) }\ 
                         A^8(F\times     F)\Bigr)\\  
                         =&       \psi_\ast        \Bigl(    A^4_{(2)}(F)\otimes A^4_{(0)}(F)\oplus A^4_{(0)}(F)\otimes A^4_{(2)}(F)\Bigr)\\
                         =&  \psi_\ast \Bigl(     \bigl( S_1 A^4(F)\otimes S_2 A^4(F)\oplus S_2 A^4(F)\otimes S_1 A^4(F)\bigr)\cap  A^8_{hom}(F\times F)\Bigr)\\
                         =&   S_3^{mob} A^8(Z)\cap A^8_{hom}(Z)\ .\\
                       \end{split}\]
   Here, the first equality is by definition of $Q$. The second equality is the fact that $P^\ast A^3_{hom}(Y)=A^2_{(2)}(F)$ (indeed, $A^2_{(2)}(F)=P^\ast P_\ast A^4_{hom}(F)$ \cite[Proof of Proposition 21.10]{SV}, and $P_\ast A^4_{hom}(F)=A^3_{hom}(Y)$, since $A^3(Y)$ is generated by lines \cite{Par}). For the third equality, one observes that $H^1(F,\QQ)=0$ and so the big line bundle $\psi^\ast(h)\in A^1(F\times F)$ can be written $\psi^\ast(h)=(p_1)^\ast(\ell_1)+(p_2)^\ast(\ell_2)$ with $\ell_j\in A^1(F)$ big and universally defined. Since the very general $F$ has Picard number 1, the $\ell_j$ must be non-zero multiples of the polarization class of $F$. Hence, one has isomorphisms $A^2_{(2)}(F)\xrightarrow{\cdot \ell_j^2} A^4_{(2)}(F)$ \cite[Proposition 22.2]{SV} and surjections 
     \[ A^2_{(0)}(F)\  \xrightarrow{\cdot \ell_j^2} \ A^4_{(0)}(F)\cong \QQ[\ell_j^4] \] (hard Lefschetz), which proves the third equality. The fourth equality is an application of \cite[Proposition 4.5]{V14}. For the last equality, one notes that 
$S_1 A^4(F)$ is generated by points lying on a uniruled divisor, and $S_2 A^4(F)$ is generated by points on the constant cycle surfaces constructed by Voisin 
\cite[Proposition 4.4]{V14}; since one can find a uniruled divisor $D$ and a constant cycle surface $S$ such that $D\times S$ is in general position with respect to the indeterminacy locus of $\psi$, this gives the last equality.   
   
   %$\psi_\ast S_j\subset S_j$ \cite[Proof of Corollary 4.9]{V14}.               
  Altogether, this proves that 
   \[ (\pi^{tr}_Z)_\ast A^8(Z)=  S_3^{mob} A^8(Z)\cap A^8_{hom}(Z)\ .\]
   
  It only remains to prove Claim \ref{Q}, i.e. we need to prove the vanishing
   \[   \Bigl( c_h\, {}^t Q\circ \Gamma_{h^6}\circ Q -\Delta_Y\Bigr)  \circ \pi^{tr}_Y=0\ \ \ \hbox{in}\ A^4(Y\times Y)\ .\]   
   Let $\pi^{prim}_Y\in A^4(Y\times Y)$ be the projector defined as
   \[ \pi^{prim}_Y:=  \Delta_Y - {1\over 3} \sum_j h^j\times h^{4-j}\ \ \   \in A^4(Y\times Y)\ .\]
   Since $H^{2,2}(Y,\QQ)$ is one-dimensional for a very general cubic $Y$, there is equality $\pi^{tr}_Y=\pi^{prim}_Y$ for very general $Y$.
   Since $(Y,\pi_Y^{tr},0)$ is a submotive of $(Y,\pi_Y^{prim},0)$, it suffices to prove the vanishing
   \[   \Bigl( c_h\, {}^t Q\circ \Gamma_{h^6}\circ Q -\Delta_Y\Bigr)  \circ \pi^{prim}_Y=0\ \ \ \hbox{in}\ A^4(Y\times Y)\ ,\]   
   for all smooth cubic fourfolds $Y$. But this cycle is universally defined (i.e., it is the restriction of a cycle in $A^4(\YY\times_B \YY)$, where $\YY\to B$ is the universal family of cubic fourfolds as in Proposition \ref{fpcubic}). Hence, thanks to the Franchetta property (Proposition \ref{fpcubic}) one is reduced to showing the vanishing
     \begin{equation}\label{vanh}  \Bigl( c_h\, {}^t Q\circ \Gamma_{h^6}\circ Q -\Delta_Y\Bigr)  \circ \pi^{prim}_Y=0\ \ \ \hbox{in}\ H^8(Y\times Y)\ ,\end{equation} 
    for the very general cubic $Y$.
    
    To prove the vanishing (\ref{vanh}), we proceed as follows:
    
    \begin{lemma}\label{tr} Let $Y\subset\PP^5$ be any smooth cubic. There exists a non-zero constant ${ c_h}$ such that
    \[   \Bigl( c_h\, {}^t Q\circ \Gamma_{h^6}\circ Q -\Delta_Y\Bigr){}_\ast=0\colon\ \ H^4_{tr}(Y)\ \to\ H^4_{tr}(Y)\ .\]
    \end{lemma}
    
    \begin{proof} Since $H^4_{tr}(Y)$ is the smallest Hodge structure for which the complexification contains $H^{3,1}(Y)\cong\C$, it suffices to prove there is a non-zero constant ${ c_h}$ such that
    \[   \Bigl( c_h\, {}^t Q\circ \Gamma_{h^6}\circ Q -\Delta_Y\Bigr){}_\ast=0\colon\ \ H^{3,1}_{}(Y)\ \to\ H^{3,1}_{}(Y)\ .\]  
    By definition of $Q$, it suffices to prove there is a non-zero constant $c_h$ making the following diagram commute:
    \[ \begin{array}[c]{ccccccc}
      H^{3,1}(Y) & \xrightarrow{P^\ast}& H^{2,0}(F) & \xrightarrow{{p_1}^\ast-{p_2}^\ast}& H^{2,0}({F\times F}) & \xrightarrow{{\psi}_\ast}& H^{2,0}(Z)\\
      &&&&&&\\
     \ \  \downarrow{\scriptstyle {1\over c_h}\ide} &&   \ \  \downarrow{\scriptstyle {\cdot (h_F)^2}}    && &&   \ \  \downarrow{\scriptstyle {\cdot h^6}}\\
     &&&&&&\\
      H^{3,1}(Y) & \xleftarrow{P_\ast}& H^{4,2}(F) & \xleftarrow{({p_1})_\ast-({p_2})_\ast}& H^{8,6}({F\times F}) & \xleftarrow{{\psi}^\ast}& H^{8,6}(Z)\\  
      \end{array}\]
      (Here $h_F\in A^1(F)$ denotes a polarization,) 
      %and $\wt{F\times F}\to F\times F$ is a resolution of indeterminacy of $\psi$, i.e. $\wt{F\times F}$ is birational to $F\times F$ and $\psi$ extends to a morphism $\wt{\psi}\colon\wt{F\times F}\to Z$. The morphisms $\wt{p_i}$ are defined by composing $\wt{F\times F}\to F\times F$ with the projections $p_i$.)
     
     The left square of this diagram is well-known to be commutative (essentially this is the Abel--Jacobi isomorphism established in \cite{BD}).
     To see that the right-hand square of this diagram commutes, we observe that thanks to Proposition \ref{psi}, we have a relation $\wt{\psi}^\ast(\omega_Z)=\wt{p_1}^\ast(\omega_F)-\wt{p_2}^\ast(\omega_F)$. This implies that the top horizontal line is an isomorphism: more precisely,
    the composition $  \psi_\ast (p_1^\ast -p_2^\ast) $ sends $\omega_F$
     to $6\omega_Z\in H^{2,0}(Z)$. 
        Let  $\omega_Z^\vee\in H^{6,8}(Z)$ be the class dual to $\omega_Z$ under cup product, and let $\bar{\omega}_Z^\vee\in H^{8,6}(Z)$ be its complex conjugate.
    By hard Lefschetz, we know that $\omega_Z\cdot h^6\in H^{8,6}(Z)$ is a rational multiple of the class $\bar{\omega}_Z^\vee$.
  By duality, the composition $((p_1)_\ast -(p_2)_\ast)  \psi^\ast$ sends $\omega_Z^\vee\in H^{6,8}(Z)$ to
      ${1\over 6}\omega_F^\vee\in H^{0,2}(F)$, and so  $\bar{\omega}_Z^\vee\in H^{6,8}(Z)$ is sent 
     to ${1\over 6}\bar{\omega}_F^\vee\in H^{0,2}(F)$. Using hard Lefschetz on $F$, we can rescale $h_F$ by some rational multiple to make the right-hand square commutative. This proves Lemma \ref{tr}.
      %By duality, it follows that the lower horizontal line is an isomorphism. The right vertical arrow is an isomorphism by hard Lefschetz. Since $H^{3,1}(Y)\cong\C$, follows that there exists a non-zero constant ${1\over c_h}\in\C$ making the diagram commute. 
%     
              
%      It remains to establish commutativity of the right square of the diagram, up to replacing $h_F$ by some non-zero multiple. The symplectic form $\omega_F$ is a generator of $H^{2,0}(F)$, and thanks to Proposition \ref{psi} we know there is a relation $\wt{\psi}^\ast(\omega_Z)=\wt{p_1}^\ast(\omega_F)-\wt{p_2}^\ast(\omega_F)$.
%     Since $\psi$ has degree $6$, it follows that
%     \[ \wt{\psi}_\ast \bigl( \wt{p_1}^\ast(\omega_F)-\wt{p_2}^\ast(\omega_F)\bigr) = 6\, \omega_Z\ \ \ \hbox{in}\ H^{2,0}(Z)\ .\]
%     This implies that
%     \[ \begin{split}  \wt{\psi}^\ast \Bigl( \wt{\psi}_\ast \bigl( \wt{p_1}^\ast(\omega_F)-\wt{p_2}^\ast(\omega_F)\bigr)\cdot h^6\Bigr) &
%          = 6        
%          \wt{\psi}^\ast(\omega_Z)\cdot \wt{\psi}^\ast(h^6)\\
%        & =     6   \bigl(  \wt{p_1}^\ast(\omega_F)-\wt{p_2}^\ast(\omega_F)\bigr)\cdot \wt{\psi}^\ast(h^6)\ .\\
%        \end{split}\]    
%        As before, we can write $\wt{\psi}^\ast(h)= \wt{p_1}^\ast(\ell_1)+\wt{p_2}^\ast(\ell_2)$, where the $\ell_i$ are non-zero and proportional to $h_F$. Thus, $\wt{\psi}^\ast(h)^6$ can be written as a polynomial $q(\wt{\ell_1}^4\wt{\ell_2}^2, \wt{\ell_1}^3\wt{\ell_2}^3, \wt{\ell_1}^2\wt{\ell_2)^4)$, where we use the shorthand
      \end{proof}
    
    Lemma \ref{tr} implies that one can write
    \[       \Bigl( c_h\, {}^t Q\circ \Gamma_{h^6}\circ Q -\Delta_Y\Bigr) = \gamma\ \ \ \hbox{in}\ H^8(Y\times Y)\ ,\]
    where $\gamma$ is a {\em completely decomposed cycle\/}, i.e. a cycle supported on a finite union of subvarieties $V_i\times W_i\subset Y\times Y$,
    with $\dim V_i+\dim W_i=4$. Since the composition of the completely decomposed cycle $\gamma$ with $\pi^{prim}_Y$ is zero (for any smooth cubic $Y$), this implies the necessary vanishing (\ref{vanh}), and so the theorem is proven.    
    \end{proof} 
  
  \begin{remark} For {\em cyclic\/} cubic fourfolds (i.e. cubic fourfolds of the form $f(x_0,\ldots,x_4)+ x_5^3=0$), a statement similar to Theorem \ref{motive} was proven in \cite[Proposition 3.7]{CCL}.
  \end{remark}

%  \begin{proposition}\label{altern} Let $Y\subset\PP^5$ be a smooth cubic not containing a plane, and let $Z:=Z(Y)$ be the associated LLSS eightfold. Let $h\in A^1(Z)$ be the polarization.
%  There is equality
%    \[ S_3 A^8(Z)\cap A^8_{hom}(Z) =. \ima\Bigl( A^2_{hom}(Z)\xrightarrow{\cdot h^6} A^8(Z)\Bigr)\ .\]  
%  \end{proposition}
%  
  
  We also define the corresponding piece in the codimension $2$ Chow group of $Z$:
  
  \begin{definition}\label{A22} For any LLSS eightfold $Z$, let  
  \[   A^2_{[2]}(Z):= Q_\ast A^3_{hom}(Y)\ \ \subset\ A^2_{hom}(Z) \]
  (where $Q$ is as in (\ref{defQ})). 
  \end{definition}
  
  The piece $ A^2_{[2]}(Z)$ is likely to coincide with $A^2_{hom}(Z)$, although this seems hard to prove (just as it is hard to prove that the inclusion $A^2_{(2)}(F(Y))\subset A^2_{hom}(F(Y))$ is an equality, cf. \cite{SV}).

 \begin{proposition}\label{ttr} Let $Y\subset\PP^5$ be a smooth cubic not containing a plane, and let $F:=F(Y)$ and $Z:=Z(Y)$ be the Fano variety of lines, resp. the associated LLSS eightfold. 
  
  \begin{enumerate}[(i)]
  
%  \item
%  One has
%   \[  \psi^\ast A^2(Z)\cap  A^2_{(2)}(F\times F) = \psi^\ast A^2_{[2]}(Z)\ .\]
  
 \item
 The projector ${}^t \pi^{tr}_Z\in A^8(Z\times Z)$ has the property that
   \[    ({}^t \pi^{tr}_Z)_\ast A^\ast(Z) =   A^2_{[2]}(Z)\ .\]
   
  \item
  There is an isomorphism
    \[ \cdot h^6\colon\ \ A^2_{[2]}(Z)\ \xrightarrow{\cong}\  S_3^{mob} A^8(Z)\cap A^8_{hom}(Z)\ \]
    (In particular, there is a correspondence-induced isomorphism $A^2_{[2]}(Z)\cong A^3_{hom}(Y)$).
    \end{enumerate}
           \end{proposition}
   
   \begin{proof}
   
%   \noindent
%   (\rom1) It follows from the proof of Claim \ref{Q} that
%     \[ A^3_{hom}(Y)\ \xrightarrow{P^\ast}\ A^2(F)\ \xrightarrow{(p_1)^\ast -(p_2)^\ast}\ A^2(F\times F) \]
%     is injective, with left-inverse given by
%     \[ A^2(F\times F)\ \xrightarrow{\psi^\ast\psi_\ast}\ A^2(F\times F)\ \xrightarrow{\cdot c_h\, h^6}\ A^8(F\times F)\ \xrightarrow{(p_1)_\ast-(p_2)_\ast}\ A^4(F)\ \xrightarrow{P_\ast}\ A^3(Y)\ .\]
%     This implies that there is equality
%     \[    (p_1^\ast-p_2^\ast)P^\ast = \psi^\ast Q_\ast\colon\ \ A^3_{hom}(Y)\ \to\ A^2(F\times F)\ .\]
%    But $ (p_1^\ast-p_2^\ast)P^\ast A^3_{hom}(Y)$ is exactly $A^2_{(2)}(F)\otimes A^0(F)\oplus A^0(F)\otimes A^2_{(2)}(F)= A^2_{(2)}(F\times F)$ (indeed, in \cite[Proof of Proposition 20.3]{SV}, it is proven that $A^2_{(2)}(F)=  P^\ast P_\ast A^4_{hom}(F)$. Since $A^3(Y)$ is generated by lines, this gives $A^2_{(2)}(F)=P^\ast A^3_{hom}(Y)$). This proves (\rom1).  
   
   \noindent
   (\rom1) Claim \ref{Q} guarantees that
     \[ {}^t \pi^{tr}_Z= c_h\, Q\circ \pi^{tr}_Y\circ {}^t Q\circ \Gamma_{h^6}\ \ \ \in A^8(Z\times Z)\]
     is idempotent, and that
     \[ Q\colon\ \ t(Y):=(Y,\pi^{tr}_Y,0)\ \to\     (Z,{}^t \pi^{tr}_Z,-1)\ \ \hbox{in}\ \MM_{\rm rat} \]
     is an isomorphism, with inverse $c_h\, {}^t Q\circ \Gamma_{h^6}$. It follows that
     \[   ({}^t \pi^{tr}_Z)_\ast A^j(Z)   =0\ \ \forall j\not=2\ ,\]
     and 
     \[  ({}^t \pi^{tr}_Z)_\ast A^2(Z)   =Q_\ast A^3_{hom}(Y)=:A^2_{[2]}(Z)\ .\]
   
   \noindent
   (\rom2) This is immediate from the proof of Theorem \ref{motive}, where it is shown that
     \[  A^3_{hom}(Y)\cong S_3^{mob} A^8(Z)\cap A^8_{hom}(Z)=   (\pi^{tr}_Z)_\ast A^8(Z)= \ima\Bigl(  Q_\ast A^3_{hom}(Y)\ \xrightarrow{\cdot h^6}\ A^8(Z)\Bigr)\ .\]   
     %Indeed, surjectivity is clear from this. As for injectivity: going back to $Y$ this also follows readily
   \end{proof}

  \subsection{The anti-symplectic involution and $0$-cycles}

  \begin{theorem}\label{main-1} Let $Z=Z(Y)$ be a LLSS eightfold, where $Y\subset\PP^5$ is any smooth cubic not containing a plane. Let 
  $\iota\in\aut(Z)$ be the anti-symplectic involution;
  Then 
    \[   \begin{split} \iota_\ast&=\ide\colon\ \ S_4 A^8_{}(Z)\ \to\ A^8(Z)\ ,\\
               \iota_\ast&=-\ide\colon\ \ S_3^{mob} A^8_{}(Z)\cap A^8_{hom}(Z)\ \to\ A^8(Z)\ .\\
               \end{split}\]
       \end{theorem}
       
       \begin{proof} Let $F=F(Y)$ be the Fano variety of lines in $Y$. Via Voisin's rational map, one sees as in \cite{V14} that $S_4A^8_{}(Z)$ is one-dimensional with generator $\psi_\ast (o_F\times o_F)$, where $o_F\in A^4_{(0)}(F)=S_2 A^4(F)$ is the distinguished zero-cycle of \cite{SV}.
       There is a commutative diagram
      \[ \begin{array}[c]{ccc}
          A^8(F\times F) &\xrightarrow{\psi_\ast}& A^8(Z)\\
          &&\\
          \downarrow{\scriptstyle( \iota_F)_\ast}&&\downarrow{\scriptstyle \iota_\ast}\\
          &&\\
         A^8(F\times F) &\xrightarrow{\psi_\ast}& \ \, A^8(Z)\ ,\\
         \end{array}\]
       where $\iota_F$ is the map exchanging the two factors  \cite[Remark 2.19]{CCL}. Since $o_F\times o_F$ is invariant under $\iota_F$, it follows that
       $\iota$ acts as the identity on $S_4 A^8_{}(Z)$.
       
       To prove the statement for $S_3^{mob} A^8(Z)\cap A^8_{hom}(Z)$, let us first ensure that $\iota_\ast$ preserves this subgroup:
       
       \begin{lemma}\label{preserve} We have
         \[ \iota_\ast  \bigl( S_3^{mob} A^8(Z)\cap A^8_{hom}(Z)\bigr)\ \ \subset\             
           S_3^{mob} A^8(Z)\cap A^8_{hom}(Z)\ .\]
           \end{lemma}

           \begin{proof} Clearly, the involution $\iota$ preserves $S_\ast$ and preserves homological triviality. To see that $\iota$ also preserves
           $S_\ast^{mob}$, let $x\in Z$ be a point outside the uniruled divisor $D$ and having a $j$-dimensional orbit. The inverse image $\psi^{-1}(x)\subset F\times F$
           consists of points outside the indeterminacy locus $I$ of $\psi$. Let $i\in\aut(F\times F)$ be the map exchanging the two factors. Since $I$ is symmetric, 
           $i_\ast(\psi^{-1}(x))$ is disjoint from $I$, and so $\iota(x)$, which is equal to the support of $\psi_\ast i_\ast(\psi^{-1}(x))$ (Remark \ref{rmk: compatibility involution})   
           lies outside of $D$, i.e. $\iota_\ast(x)\in S_j^{mob} A^8(Z)$.
                   %           As we have seen, 
%           \[ A^8_{(2)}(Z)= A^2_{hom}(Z)\cdot h^6\ ,\]
%           where $h$ is any ample divisor (\ref{}). Taking $h$ a $\iota$-invariant ample divisor, this proves the Lemma.
            \end{proof}
           
           Next, we use the fact that
           \[  \begin{split} S_3^{mob}\cap A^8_{hom}(Z)&=    (\Gamma_{h^6}\circ Q)_\ast A^3_{hom}(Y)\\        
                                          &=  (\Gamma_{h^6}\circ Q\circ \pi^{tr}_Y)_\ast A^3_{}(Y)\ ,\\     
                                          \end{split}\]    
                              (where $Q$ is the correspondence of (\ref{defQ}), and $\pi^{tr}_Y$ is the projector on $H^4_{tr}(Y)$ as in Proposition \ref{bp}). The involution $\iota$, being anti-symplectic, acts as $-\ide$
                              on $H^2_{tr}(Z) =  (Q\circ \pi^{tr}_Y)_\ast H^4(Y)$, and so
                        \[      c_h\,\,    {}^t Q\circ \Gamma_\iota\circ\Gamma_{h^6}\circ   Q\circ \pi^{tr}_Y  + \pi^{tr}_Y =0\ \ \ \hbox{in}\  H^8(Y\times Y)      \]
                        (where $c_h\in\QQ^\ast$ and $h$ are as in Claim \ref{Q}).
                        
   Let us now consider things family-wise. As in Proposition \ref{fpcubic}, let $\YY\to B$ and $\Zz\to B$ denote the universal family of smooth cubic fourfolds not containing a plane, resp. LLSS eightfolds. Let $Q\in A^3(\YY\times_B \Zz)$ be the relative correspondence as in the proof of Claim \ref{Q}, and let $\Gamma_{h^6}\in A^{14}(\Zz\times_B \Zz)$ be the relative correspondence induced by some ample element $h\in A^1(\Zz)$. Let $h_\YY\in A^1(\YY)$ be a relative hyperplane section. The correspondence
   \[ \pi^{prim}_\YY:=\Delta_\YY - {1\over 3}\sum_{j=0}^4  (pr_1)^\ast (h_\YY)^j\cdot   (pr_2)^\ast (h_\YY)^{4-j}\ \ \ \in\ A^4(\YY\times_B \YY)   \]
   is such that $\pi^{prim}_\YY\vert_b=\pi^{tr}_{Y_b}$ for $b\in B$ very general, and
   \begin{equation}\label{prim} (\pi^{prim}_\YY\vert_b )_\ast= (\pi^{tr}_{Y_b})_\ast=\ide\colon\ \ A^3_{hom}(Y_b)\ \to\ A^3(Y_b)\ \ \ \forall b\in B\ .\end{equation}
   It follows from the above considerations that the relative correspondence   
   \[ c_h\,\,    {}^t Q\circ \Gamma_\iota\circ\Gamma_{h^6}\circ   Q\circ \pi^{prim}_\YY  + \pi^{prim}_\YY \ \ \in\ A^4(\YY\times_B \YY)   \]
   is homologically trivial when restricted to the very general fibre. But then, the Franchetta property (Proposition \ref{fpcubic}) ensures that this relative correspondence is fibrewise zero, i.e.
   \[    \Bigl(c_h\,\,    {}^t Q\circ \Gamma_\iota\circ\Gamma_{h^6}\circ   Q\circ \pi^{prim}_\YY\Bigr)\vert_b  = - \pi^{prim}_\YY\vert_b \ \ \in\ A^4(Y_b\times Y_b) \ \ \ \forall b\in B\ .\]
   Composing on both sides with $Q\vert_b$ on the right and with $(c_h\, \Gamma_{h^6}\circ Q\circ\pi^{prim}_\YY)\vert_b$ on the left, it follows that
   \[ (c_h)^2\Bigl(\Gamma_{h^6}\circ Q\circ  \pi^{prim}_\YY\circ  {}^t Q\circ \Gamma_\iota\circ\Gamma_{h^6}\circ   Q\circ \pi^{prim}_\YY\circ {}^t Q\Bigr)\vert_b    =   -c_h\Bigl( \Gamma_{h^6}\circ Q\circ \pi^{prim}_\YY\circ{}^t Q\Bigr)\vert_b\ \ \ \in A^8(Z_b\times Z_b) \ ,\]
   for all $b\in B$.
   Applying this equality of correspondences to $A^8_{hom}(Z_b)$, and using (\ref{prim}), we obtain the equality of actions
   \[  (\pi^{tr}_{Z_b}\circ \Gamma_{\iota_b}\circ \pi^{tr}_{Z_b})_\ast = -(\pi^{tr}_{Z_b})_\ast\colon\ \ A^8_{hom}(Z_b)\ \to\ A^8(Z_b)\ \ \ \forall b\in B\ ,\]
   where $\pi^{tr}_{Z_b}:=c_h\,\Gamma_{h_b^6}\circ Q_b\circ \pi^{tr}_{Y_b}\circ{}^t Q_b$ is the projector on $S_3^{mob}\cap A^8_{hom}(Z_b)$ as in Theorem
   \ref{motive}.
   Combined with Lemma \ref{preserve}, this means that
   \[ (\iota_b)_\ast=-\ide\colon\ \ S^{mob}_3 A^8_{}(Z_b)\cap A^8_{hom}(Z_b)\ \to\ A^8(Z_b)\   \ \ \forall b\in B\ ,\]
   and so the theorem is proven.
          \end{proof}

  \begin{theorem}\label{main0} Let $Z$ be a LLSS eightfold birational to a Hilbert scheme $Z^\prime=S^{[4]}$ where $S$ is a K3 surface.  
  Let $\iota^\prime\in\bir(Z^\prime)$ be the map induced by $\iota\in\aut(Z)$. Then
    \[  \begin{split} (\iota^\prime)_\ast&=\ide\colon\ \ A^8_{(j)}(Z^\prime)\ \to\ A^8_{}(Z^\prime)\ \ \ \hbox{for}\ j\in\{0,4,8\}\ ,\\
               (\iota^\prime)_\ast &=-\ide\colon\ \ \   A^8_{(j)}(Z^\prime)  \ \to\ A^8_{}(Z^\prime)   \ \ \ \hbox{for}\ j\in\{2,6\}\ .\\
               \end{split}\]   
        (Here the bigrading refers to the MCK decomposition of \cite{V6}.)       
       \end{theorem}
  
  \begin{proof} The Hilbert scheme $Z^\prime$ has a MCK decomposition \cite{V6}, and the induced decomposition of $A^8(Z^\prime)$ is compatible with Voisin's orbit filtration $S_\ast$, in the sense that
    \begin{equation}\label{inc} \bigoplus_{j=0}^r A^8_{(2j)}(Z^\prime) = S_{4-r}  A^8(Z^\prime) \ \end{equation}
  \cite[End of section 4.1]{V14}. 
  
 Using Rie\ss's work (Theorem \ref{rie}), the MCK decomposition of $Z^\prime$ induces a MCK decomposition for $Z$, and there exists a correspondence $R_\phi\in A^8(Z^\prime\times Z)$ inducing an isomorphism of bigraded rings $A^\ast_{(\ast)}(Z^\prime)\cong A^\ast_{(\ast)}(Z)$.  
In order to prove Theorem \ref{main0}, we want to relate the orbit filtration $S_\ast$ on $A^8(Z)$ with the ``motivic'' decomposition $A^8_{(\ast)}(Z)$. This is the content of the following lemma.

%  
%  As we have seen, there is equality
%   \begin{equation}\label{this} S_3\cap A^8_{hom}(Z)= A^2_{[2]}(Z)\cdot h^6\ ,\end{equation}
%           where $h$ is the polarization (Proposition \ref{ttr}). 

%           On the other hand, we also have equality
%           \[    A^8_{(2)}(Z^\prime) = A^2_{(2)}(Z^\prime)\cdot (h^\prime)^6\ ,\]
%           for any big divisor $h^\prime\in A^1(Z^\prime)$, and hence in particular
%           \[ \begin{split}A^8_{(2)}(Z^\prime) &= A^2_{(2)}(Z^\prime)\cdot \bigl((R_\phi)^\ast(h)\bigr)^6\ \\
%&= A^2_{hom}(Z^\prime)\cdot \bigl((R_\phi)^\ast(h)\bigr)^6\ ,\\
% \end{split}\]
%          since 
%          $(R_\phi)^\ast(h)$ is simply the birational transform of $h$ and so it is a big divisor on $Z^\prime$.   
%          Going back to $Z$ via Rie\ss's correspondence $R_\phi$, this implies that
%          \begin{equation}\label{this2}  A^8_{(2)}(Z)=(R_\phi)_\ast  A^8_{(2)}(Z^\prime) =      A^2_{(2)}(Z)\cdot     h^6\ .\end{equation}     
%           Comparing (\ref{this}) and (\ref{this2}) (and remembering the inclusion $A^2_{(2)}(Z)\subset A^2_{hom}(Z)$, which conjecturally is an equality), we find that there is an inclusion (which conjecturally is an equality)

\begin{lemma}\label{incl} Let $Z$ be as in Theorem \ref{main0}, and let $A^\ast_{(\ast)}(Z)$ refer to the bigrading induced by the bigrading on $A^\ast(Z^\prime)$. There is an inclusion
           \[ A^8_{(2)}(Z) \ \subset\  S_3^{mob}\cap A^8_{hom}(Z)\ .\]   
       \end{lemma}    
       
      \begin{proof} Let $\phi\colon Z^\prime\dashrightarrow Z$ denote the birational map. By virtue of Rie\ss's work (Theorem \ref{rie}), there is equality
        \[  A^8_{(2)}(Z)= (\bar{\Gamma}_\phi)_\ast A^8_{(2)}(Z^\prime)\ .\]
        Thanks to equality (\ref{inc}), it follows there is an inclusion
        \begin{equation}\label{aninc} A^8_{(2)}(Z)\ \subset\  (\bar{\Gamma}_\phi)_\ast  S_3 A^8(Z^\prime)\ .\end{equation}
        Let $I\subset Z^\prime$ denote the indeterminacy locus of the map $\phi$, and let $D^\prime\subset Z^\prime$ denote the strict transform of the uniruled divisor $D\subset Z$ coming from the rational map $\psi$. For any point $z\in Z^\prime$, the maximal-dimensional components of the orbit under rational equivalence $O_z$ are Zariski dense in $Z^\prime$ (\cite[Lemma 3.5]{V12}, cf. also \cite[Proof of Theorem 0.5]{SYZ} where this precise statement is obtained). 
        %(Actually, one could use Thm 2.6 of Voisin's "Remarks on coisotropic" to see that a point in $S_3$ is supported on $C^{(3)}$ for some CC curve $C$...)
        Hence, the orbit of any point in $Z^\prime$ has a maximal-dimensional irreducible component not contained in $I\cup D^\prime$, and hence
        \begin{equation}\label{aninc2}    (\bar{\Gamma}_\phi)_\ast  S_3 A^8(Z^\prime)\ \subset\  S_3^{mob} A^8(Z)  \ .\end{equation}
        Combining (\ref{aninc}) and (\ref{aninc2}), one obtains
        \[   A^8_{(2)}(Z)\ \subset\  S_3^{mob} A^8(Z)  \ .\]                
        Since clearly $A^8_{(2)}(Z)\subset A^8_{hom}(Z)$, this proves the lemma.      
        \end{proof} 
           
         %  as it should. 
Theorem \ref{main-1}, combined with Lemma \ref{incl}, implies that 
             \begin{equation}\label{ONZ} \iota_\ast =-\ide\colon\ \ A^8_{(2)}(Z)\ \to A^8(Z)\ .\end{equation}
%      Taking a $\iota$-invariant ample divisor and invoking the isomorphism $\cdot h^6\colon A^2_{(2)}(Z)\cong   A^8_{(2)}(Z)$, it follows that also
%       \begin{equation}\label{a2} \iota_\ast =-\ide\colon\ \ A^2_{(2)}(Z)\ \to A^2(Z)\ .\end{equation}    
%          

%            Because $Z^\prime$ is a Hilbert scheme, the intersection product induces surjections
%   \[  A^2_{(2)}(Z^\prime)^{\otimes  j}\ \twoheadrightarrow\ A^{2j}_{(2j)}(Z^\prime)\ \ \ (j=1,2, 3, 4)\ \]
%   (Proposition \ref{mult}). Applying the isomorphism $(R_\phi)_\ast$, it follows there are also surjections
%    \[  A^2_{(2)}(Z)^{\otimes  j}\ \twoheadrightarrow\ A^{2j}_{(2j})(Z)\ \ \ (j=1,2, 3, 4)\ .\]
%    Because intersection product and $\iota^\ast$ commute, equality (\ref{a2}) now implies that
%     \[ \iota_\ast =(-1)^j\, \ide\colon\ \ A^{2j}_{(2j)}(Z)\ \to A^{2j}(Z)\ \ \ (j=1,2, 3, 4)\ .\]
%     Again intersecting with a $\iota$-invariant ample divisor, this gives the equality of actions on $0$-cycles
%      \begin{equation}\label{a3} \iota_\ast =(-1)^j\, \ide\colon\ \ A^{8}_{(2j)}(Z)\ \to A^{8}(Z)\ \ \ (j=1, 2, 3, 4)\ .\end{equation}  
      
    There is a commutative diagram
    \begin{equation}\label{comdi} \begin{array}[c]{ccc}
       A^8(Z^\prime) &\xrightarrow{(R_\phi)_\ast}& A^8(Z)\\
       &&\\
      \ \  \downarrow{\scriptstyle (\iota^\prime)_\ast}      && \ \  \downarrow{\scriptstyle \iota_\ast}\\
      &&\\
       A^8(Z^\prime) &\xrightarrow{(R_\phi)_\ast}& A^8(Z)\\
       \end{array}\end{equation}
       (To check commutativity of this diagram, we note that (thanks to the moving lemma) any $0$-cycle $b\in A^8(Z^\prime)$ can be represented by a cycle $\beta$ supported on an open $U^\prime$ such that the birational map $\phi\colon Z^\prime\dashrightarrow Z$ restricts to an isomorphism $\phi\vert_{U^\prime}\colon U^\prime\cong U$ and $U\subset Z$ is $\iota$-stable. As $(R_\phi)_\ast$ and 
 $(\bar{\Gamma}_\phi)_\ast$ coincide on $0$-cycles (Theorem \ref{rie}), the image $(R_\phi)_\ast(b)\in A^8(Z)$ is represented by the cycle $(\phi\vert_{U^\prime})_\ast (\beta)$ with support on $U$. This shows commutativity of the diagram.)
 
 This commutative diagram, plus the fact that $\iota$ acts as the identity on $S_4 A^8(Z)=A^8_{(0)}(Z)$ (Theorem \ref{main-1}), proves the $j=0$ case of the theorem. Moreover, the commutative diagram, plus equality (\ref{ONZ}), gives 
   \begin{equation}\label{ONP} 
       ( \iota^\prime)_\ast =-\ide\colon\ \ A^8_{(2)}(Z^\prime)\ \to A^8(Z^\prime)\ ,\end{equation}      
     which is the $j=2$ case of the theorem.   
     
     Taking a $\iota$-invariant ample divisor $h\in A^1(Z)$ and invoking the isomorphism 
     \begin{equation}\label{hlt} \cdot h^6\colon \ \ A^2_{(2)}(Z)\cong   A^8_{(2)}(Z)\end{equation}
     (Remark \ref{birat2}), it follows from (\ref{ONZ}) that also
       \begin{equation}\label{a2} \iota_\ast =-\ide\colon\ \ A^2_{(2)}(Z)\ \to A^2(Z)\ .\end{equation}    
  The intersection product induces surjections
   \[  A^2_{(2)}(Z)^{\otimes  j}\ \twoheadrightarrow\ A^{2j}_{(2j)}(Z)\ \ \ (j=1,2, 3, 4)\ \]
   (Remark \ref{birat2}).
   Since intersection product and $\iota^\ast$ commute, equality (\ref{a2}) thus implies that 
     \[ \iota_\ast =(-1)^j\, \ide\colon\ \ A^{2j}_{(2j)}(Z)\ \to A^{2j}(Z)\ \ \ (j=1,2, 3, 4)\ .\]
     Taking once again a $\iota$-invariant ample divisor $h\in A^1(Z)$, and invoking the isomorphism 
     \[ \cdot h^{8-2j}\colon \ \ A^{2j}_{(2j)}(Z)\cong   A^8_{(2j)}(Z)\] 
     (Remark \ref{birat2}), it follows that
  \begin{equation}\label{almost} \iota_\ast =(-1)^j\, \ide\colon\ \ A^{8}_{(2j)}(Z)\ \to A^{8}(Z)\ \ \ (j=1,2, 3, 4)\ .\end{equation}
 Combining (\ref{almost}) with the commutative diagram (\ref{comdi}), we find that
   \[  ( \iota^\prime)_\ast =(-1)^j\, \ide\colon\ \ A^{8}_{(2j)}(Z^\prime)\ \to A^{8}(Z^\prime)\ \ \ (j=1, 2, 3, 4)\ .  \]
  The theorem is now proven.  
    \end{proof}

  \subsection{A reformulation} We can reformulate Theorem \ref{main0} into a simpler statement, that does not mention MCK decompositions or bigradings on the Chow ring:
  
  \begin{theorem}\label{main2}  Let $Z$ be a LLSS eightfold birational to a Hilbert scheme $Z^\prime=S^{[4]}$ where $S$ is a K3 surface.  
  Let $\iota^\prime\in\bir(Z^\prime)$ be the map induced by $\iota\in\aut(Z)$. Then
         \[ \begin{split} (\iota^\prime)_\ast [x,y,z,t]= [x,y,z,t] &-2\bigl( [x,y,z,o] + [x,y,t,o]+[x,z,t,o] + [y,z,t,o]\bigr)\\ &+4\bigl( [x,y,o,o]+[x,z,o,o]+\ldots+[z,t,o,o]\bigr)\\ &-8\bigl( [x,o,o,o]+[y,o,o,o]+[z,o,o,o]+[t,o,o,o]\bigr)\\&+16[o,o,o,o]\ .\\
     \end{split}\]
     (Here, $x,y,z,t\in S$ and $o\in A^2(S)$ denotes the distinguished $0$-cycle.)
\end{theorem}

\begin{proof} This is a translation of Theorem \ref{main0}, using the generators for the various pieces $A^8_{(j)}(Z^\prime)$ given in Theorem \ref{K3m}.
We know that $\iota^\prime$ preserves these pieces (Theorem \ref{main0}), and so one can check the formula of Theorem \ref{main2} piece by piece. This is a direct computation. For example, Theorem \ref{K3m} says that $A^8_{(8)}(Z^\prime)$ is generated by expressions
of the form
  \[ \begin{split} [x,y,z,t] - \bigl( [x,y,z,o] + [x,y,t,o] + \cdots\bigr) + \bigl( [x,y,o,o] + [x,z,o,o]+\cdots\bigr)&\\ - \bigl( [x,o,o,o] + [y,o,o,o]+\cdots\bigr) + [o,o,o,o]&\ ,\\
  \end{split}\]
  and the formula of Theorem \ref{main2} leaves expressions of this form invariant, which is in agreement with the fact that $\iota^\prime$ is the identity on $A^8_{(8)}(Z^\prime)$ (Theorem \ref{main0}). The verification of the action of $\iota^\prime$ on $A^8_{(j)}(Z^\prime)$, $j<8$ is similar (but easier).
\end{proof}

 \subsection{Another reformulation}
 
 Theorem \ref{main0} also admits a reformulation that does not mention LLSS eightfolds:
 
 \begin{theorem}\label{main25} Let $S$ be a K3 surface of Picard number 1 and degree $d$ satisfying $d>2$ and $d= (6n^2-6n+2)/a^2$ for some $n,a\in\ZZ$. Let $X=S^{[4]}$ be the Hilbert scheme. Then $\bir(X)=\ZZ/2\ZZ$ and the non-trivial element $\iota\in\bir(X)$ acts on $A^8(X)$ as in Theorem \ref{main2}.
 \end{theorem}
 
 \begin{proof} Given a K3 surface $S$ of Picard number 1 and degree $d>2$, the Hilbert scheme $X=S^{[4]}$ has $\bir(X)$ a group of order $1$ or $2$
 \cite[Section 4.3.1]{debarre}.
 On the other hand, the condition $d= (6n^2-6n+2)/a^2$ for some $n,a\in\ZZ$ implies (and actually is equivalent to the fact) that $X$ is birational to an LLSS eightfold \cite[Proposition 5.3]{LPZ}, and so $\bir(X)$ is non-trivial.
 It follows that $X$ as in the theorem must have $\bir(X)=\ZZ/2\ZZ$, and the non-trivial element 
 $\iota\in\bir(X)$ acts on $A^8(X)$ as in Theorem \ref{main2}.
  \end{proof}

  \subsection{A consequence}
  
  \begin{corollary}\label{cor} Let $Z$ be a LLSS eightfold birational to a Hilbert scheme $Z^\prime=S^{[4]}$ where $S$ is a K3 surface. Let $X:=Z/\langle\iota\rangle$ be the quotient under the anti-symplectic involution $\iota\in\aut(Z)$.
Then
  \[ \begin{split} &\ima \Bigl(  A^2(X)^{\otimes 4}\to A^8(X)\Bigr)
                          = \ima\Bigl( A^3(X)\otimes A^2(X)\otimes A^2(X)\otimes A^1(X)\Bigr)
                          =\QQ[h^8]\ ,\\
                          \end{split}\]
        where $h\in A^1(X)$ is an ample divisor.
    \end{corollary}
    
    \begin{proof} Recall that $X$ is a quotient variety, and so the Chow groups with $\QQ$-coefficients of $X$ have a ring structure (cf. section \ref{ssquot}).
    
   As noted before, the birationality $Z\to Z^\prime$ induces an isomorphism of bigraded rings $A^\ast_{(\ast)}(Z)\cong A^\ast_{(\ast)}(Z^\prime)$ (Theorem \ref{rie}).
%   It follows that there is a ``hard Lefschetz type'' isomorphism
%  \begin{equation}\label{hard}  \cdot h^6\colon\ \ \ A^2_{(2)}(Z)\ \xrightarrow{\cong}\ A^8_{(2)}(Z)\ ,\end{equation}
%  for any ample divisor class $h\in A^1(Z)$ (cf. Proposition \ref{mult}). Taking a $\iota$-invariant ample class $h$, this gives in particular an isomorphism
%  \[  \cdot h^6\colon\ \ \ A^2_{(2)}(Z)^{ \iota}\ \xrightarrow{\cong}\ A^8_{(2)}(Z)^{ \iota }\ ,\] 
%  and so Theorem \ref{main0} implies that
We have seen (equality (\ref{a2})) that $\iota$ acts as $-\ide$ on $A^2_{(2)}(Z)$, and so
  \[  A^2_{(2)}(Z)^{ \iota }=0\ .\]
  Similar reasoning also shows that
   \[ A^2(Z)^{ \iota }\ \subset\ A^2_{(0)}(Z)\ .\]
   Indeed, let $b\in A^2(Z)^{\iota }$, and let $h\in A^1(Z)$ be once again a $g$-invariant ample class. Then $b\cdot h^6\in A^8(Z)^{ \iota }$ and so (in view of Theorem \ref{main0}) we know that $b\cdot h^6\in A^8_{(0)}(Z)$. Writing the decomposition $b=b_0+b_2$, where $b_j\in A^2_{(j)}(Z)$, and invoking the hard Lefschetz type isomorphism (\ref{hlt}), we see that $b_2=0$.
   
   The corollary is now readily proven. For instance, the image of the intersection map
   \[ \ima \Bigl( A^3(Z)^{ \iota }\otimes A^2(Z)^{ \iota }\otimes A^2(Z)^{ \iota }\otimes A^1(Z)^{\iota }\to\ A^8(Z)\Bigr) \]   
   is contained in
   \[ \ima  \Bigl( \bigl( \bigoplus_{j\le 2} A^3_{(j)}(Z)\bigr)\otimes A^2_{(0)}(Z)\otimes A^2_{(0)}(Z)\otimes A^1(Z)\ \to\ A^8(Z)\Bigr) \ \ \subset\  \bigoplus_{j\le 2} A^8_{(j)}(Z)\ .\]
   On the other hand, the image is obviously $\iota $-invariant, and so
   \[     \ima \Bigl( A^3(Z)^{ \iota }\otimes A^2(Z)^{ \iota }\otimes A^2(Z)^{ \iota }\otimes A^1(Z)^{\iota }    \to\ A^8(Z)\Bigr)\ \ \subset\ A^8(Z)^{\iota} \subset\ A^8_{(0)}(Z)\oplus A^8_{(4)}(Z) \oplus A^8_{(8)}(Z)\]   
   (where the second inclusion is equality (\ref{almost})). It follows that
     \[   \ima \Bigl( A^3(Z)^{ \iota }\otimes A^2(Z)^{ \iota }\otimes A^2(Z)^{ \iota }\otimes A^1(Z)^{\iota }    \to\ A^8(Z)\Bigr)  \ \ \subset\ A^8_{(0)}(Z)\cong\QQ[h^8]\ . \]   
      \end{proof}

 \vskip1cm
\begin{nonumberingt} Thanks to the referee for many insightful comments. Thanks to Christian Lehn, who very patiently explained details of the involution $\iota$ on a crowded train \cite{Le}. I am also very grateful to Chiara Camere, Alberto Cattaneo, Lie Fu, Laura Pertusi, Ulrike Rie\ss\  and Charles Vial for various inspiring exchanges.

\end{nonumberingt}

\end{document}